[1994/06/01]
\documentclass[12pt]{amsart}

        [1994/05/12 v2.2b Standard AMS font definitions]
\DeclareFontFamily{U}{msa}{}
\DeclareFontShape{U}{msa}{m}{n}
  { <5> <6> <7> <8> <9> gen * msam
    <10> <10.95> <12> <14.4> <17.28> <20.74> <24.88> msam10}{}
        [1994/05/12 v2.2b Standard AMS font definitions]
\DeclareFontFamily{U}{msb}{}
\DeclareFontShape{U}{msb}{m}{n}
  { <5> <6> <7> <8> <9> gen * msbm
    <10> <10.95> <12> <14.4> <17.28> <20.74> <24.88> msbm10}{}

\newtheorem{cor}[subsection]{Corollary}

\newtheorem{prop}[subsection]{Proposition}
\newtheorem{rem}[subsection]{Remark}
\newtheorem{thm}[subsection]{Theorem}

\newtheorem{firstthm}{Theorem 1}

\newtheorem{secondthm}{Theorem 2}

\theoremstyle{definition}

\numberwithin{equation}{section}

\title[Automorphisms of Surface Braid Groups]{Automorphisms of Surface Braid Groups}
\author[E. Irmak]{Elmas Irmak}
\address{Department of Mathematics \ University of Michigan \ Ann Arbor, MI 48109-1109}
\email{eirmak@umich.edu}

\author[N. V. Ivanov]{Nikolai V. Ivanov}
\address{Department of Mathematics, Michigan State University, East Lansing, MI
48824-1027}
\email{ivanov@@math.msu.edu}
\thanks{The first author was partially supported by NSF Grant DMS 9401284.}

\author[J. D. McCarthy]{John D. McCarthy}
\address{Department of Mathematics, Michigan State University, East Lansing, MI
48824-1027}
\email{mccarthy@@math.msu.edu}
\thanks{The second author was partially supported by NSF Grant DMS 9305067.}

\date{June 3, 2003}


\newcommand{\abstracttext}{In this paper, we prove that each automorphism of a surface braid group is induced by a homeomorphism of the underlying surface,
provided that this surface is a closed, connected, orientable surface of genus at least 2, and the number of strings is at least three. This result generalizes
previous results for classical braid groups, mapping class groups, and Torelli groups.}

\begin{document}

\maketitle

\begin{abstract}  \abstracttext  \end{abstract}

\section{Introduction}
\label{sec:introduction}

Let $S$ be a compact orientable surface. The {\it Teichm\"{u}ller modular group} $Mod_S$ of $S$, also known as the {\it mapping class group} of $S$, is
the group of isotopy classes of orientation-preserving diffeomorphisms $S \rightarrow S$. The {\it pure modular group} $PMod_S$ is the subgroup of
$Mod_R$ consisting of isotopy classes of orientation-preserving diffeomorphisms $S \rightarrow S$ which preserve each component of the boundary
$\partial S$ of $S$. The {\it extended modular group} $Mod*(S)$ of $S$ is the group of isotopy classes of {\it all} (including
orientation-reversing) diffeomorphisms $S \rightarrow S$. 

Let $n$ be a positive integer. Let $R$ be the surface obtained from $S$ by removing $n$ distinct points from $S$. Note that there is a natural
homomorphism $\mu: PMod_R \rightarrow PMod_S$ corresponding to extending homeomorphisms $F : R \rightarrow R$ to homeomorphisms 
$\overline{F} : S \rightarrow S$. We denote the kernel of $\mu$ as $Mod_R(S)$ and its intersection with $PMod_R$ as $PMod_R(S)$. 

Let $S^{[n]}$ denote the space of $n$-tuples $(x_1,...,x_n)$ of distinct points of $S$. Note that the symmetric group $\Sigma_n$ acts on $S^{[n]}$ by
permuting coordinates. We recall that the $n$-string braid group $B_n(S)$ of $S$ is the fundamental group of the quotient of $S^{[n]}$ by this action of
$\Sigma_n$ (\cite{i3}). Likewise, the pure $n$-string braid group $PB_n(S)$ of $S$ is defined to be the fundamental group of $S^{[n]}$. 

Suppose $S$ has negative euler characteristic. Then the groups $Mod_R(S)$ and $PMod_R(S)$ are naturally isomorphic to $B_n(S)$ and $PB_n(S)$. The natural
isomorphism from $PB_n(S)$ to $PMod_R(S)$ arises from the connecting homomorphism in the homotopy long exact sequence of the fibration $ev : Diff(S)
\rightarrow S^[n]$ defined by the rule $ev(f) = (f(x_1),....,f(x_n))$, where $(x_1,...,x_n)$ is the chosen base point for the fundamental group $PB_n(S)$
of $S^{[n]}$. The proof that the connecting  homomorphism is injective uses the contractibility of the identity component of $Diff(S)$ (\cite{ee}). 

Hence, when $S$ has negative euler characteristic, we shall refer to $Mod_R(S)$ and $PMod_R(S)$ as the $n$-string braid group on $S$ and the pure
$n$-string braid group on $S$. If $S$ has nonnegative euler characteristic, then $Mod_R(S)$ and $PMod_R(S)$ are still closely related to the $n$-string
braid group and the pure $n$-string braid group on $S$. For instance, if $S$ is a sphere and $n \geq 3$, then $Mod_R(S)$ is naturally isomorphic to the
quotient of the $n$-string braid group $B_n(S)$ by its center, $ZB_n(S)$, which is cyclic of order two \cite{b}.

This paper concerns automorphisms of $Mod_R(S)$ and $PMod_R(S)$. The two main results of this paper generalize previous results for
automorphisms of classical braid groups, mapping class groups, and Torelli groups (\cite{dg}, \cite{i1}, \cite{m}, \cite{f}, \cite{mv}). 

\begin{firstthm}  Let $S$ be a closed connected orientable surface of positive genus. Let $n$ be an integer greater than $2$.
Let $R$ be a surface obtained from $S$ by removing $n$ distinct points from $S$. Let $\chi : PMod_R(S) \rightarrow PMod_R(S)$ be an
automorphism of the pure braid group $PMod_R(S)$. Then there exists a homeomorphism $K: R \rightarrow R$ of $R$ such that $\chi : PMod_R(S)
\rightarrow PMod_R(S)$ is equal to the automorphism $K_* : PMod_R(S) \rightarrow PMod_R(S)$ which is defined by the rule $K_*([H]) = [K
\circ H \circ K^{-1}]$ for each mapping class $[H]$ in $PMod_R(S)$. 
\end{firstthm}

\begin{secondthm} Let $S$ be a closed connected orientable surface of genus greater than $1$. Let $n$ be an integer greater than $2$.
Let $R$ be a surface obtained from $S$ by removing $n$ distinct points from $S$. Let $\chi : Mod_R(S) \rightarrow Mod_R(S)$ be an automorphism 
of the braid group $Mod_R(S)$. Then there exists a homeomorphism $F: R \rightarrow R$ of $R$ such that $\chi : Mod_R(S) \rightarrow
Mod_R(S)$ is equal to the automorphism $F_* : Mod_R(S) \rightarrow Mod_R(S)$ which is defined by the rule $F_*([H]) = [F \circ H \circ
F^{-1}]$ for each mapping class $[H]$ in $Mod_R(S)$. 
\end{secondthm} 

The authors plan to consider the case of surfaces with boundary in a future paper. 

Theorems 1 and 2 hold for almost all $n$ for all closed connected orientable surfaces $S$ regardless of
the assumption on positivity of genus. The situations not covered by the arguments used in this paper to prove Theorems $1$
and $2$ can be handled by rather straightforward arguments from known results. 

Suppose, for instance, that either (i) the genus of $S$ is equal to $0$ and $n < 2$ or (ii) the genus of $S$ is arbitary and $n = 0$ or (ii)
the genus of $S$ is equal to $1$ and $n = 1$. Then, the braid group $Mod_R(S)$ and the pure braid group $PMod_R(S)$ are both trivial. Hence,
Theorems $1$ and $2$ hold, for rather trivial reasons, when either (i), (ii), or (iii) hold. 

Suppose that the genus of $S$ is equal to $0$ and $n = 2$. Then $Mod_R(S)$ is a cyclic group of order $2$ and $PMod_R(S)$ is trivial.
Hence, every automorphism of $Mod_R(S)$ and $PMod_R(S)$ is trivial. Hence, again, Theorems $1$ and $2$ hold,
for rather trivial reasons, in this situation.

Suppose that the genus of $S$ is equal to $0$ and $n = 3$. Then $Mod_R(S)$ is isomorphic to the symmetric group $\Sigma_3$ and $PMod_R(S)$
is trivial. Hence, it is easy to see that Theorems $1$ and $2$ also hold in this situation.

Suppose that the genus of $S$ is equal to $0$. Then $Mod_S$ is trivial and hence, $Mod_R(S) = Mod_R$. It follows that $Mod_R(S)$
and $PMod_R(S)$ are both subgroups of finite index in $Mod_R$. Hence, by Theorem $3$ of \cite{k}, Theorems $1$ and
$2$ hold for all integers $n \geq 5$. 

Suppose that the genus of $S$ is at least $2$ and $n = 1$. Then $PMod_R(S) = Mod_R(S)$ and $Mod_R(S)$ is  naturally isomorphic to the
fundamental group $\pi_1(S, x)$, where $x$ is the unique puncture of $R$ on $S$. It is a well-known result that automorphisms of fundamental
groups of closed surfaces are induced by homeomorphisms (cf. \cite{zvc}, Theorem 3.3.11). Hence, Theorems $1$ and
$2$ also hold, when the genus of $S$ is at least $2$ and $n = 1$.

Hence, the only possible exceptions to Theorem $1$, for closed connected orientable surfaces $S$,
are when either (i) the genus of $S$ is equal to $0$ and $n = 4$ or (ii) the genus of $S$ is positive and $n = 2$. 
Likewise, the only possible exceptions to Theorem $2$, for closed connected orientable surfaces $S$,
are when either (i) the genus of $S$ is equal to $0$ and $n = 4$, (ii) the genus of $S$ is positive and $n = 2$, or (iii) 
the genus of $S$ is equal to $1$ and $n \geq 3$.

Theorem $2$ follows by a fairly elementary argument from Theorem $1$ using Ivanov's result that
$PMod_R(S)$ is a characteristic subgroup of $Mod_R(S)$ (\cite{i3}). This argument is given in Section \ref{sec:mainresults}. 

The bulk of the paper is devoted to proving Theorem $1$. As with the previous results mentioned above, the technique for
proving this result is to characterize algebraically certain natural elements inside $PMod_R(S)$. In the case of the foundational
result for the mapping class group $Mod_S$ of \cite{i1}, these elements were the Dehn twists in $Mod_S$. In the case of Farb's result for
the Torelli group $Tor_S$ of closed surfaces of genus greater than $3$ (\cite{f}) and the generalization of this result by McCarthy-Vautaw
to closed surfaces of genus greater than $2$ (\cite{mv}), these elements were the Dehn twists about separating circles and bounding pair
maps in $Tor_S$. 

In the present paper, the elements in question are those which we refer to as {\it $2$-string twist braids} and {\it
$1$-string bounding pair braids}. These elements constitute two of the three types of {\it basic braids} in $PMod_R(S)$ algebraically characterized
by the property that the center of their centralizers in $PMod_R(S)$ have rank equal to one, {\it $k$-string twist braids}, {\it
$k$-string bounding pair braids}, and {\it basic single-$pA$ braids}. We characterize algebraically each of these three types
of basic braids. In addition, we characterize algebraically the number of strings $k$ in twist braids and bounding pair braids.

Of particular interest to us, are the $1$-string bounding pair braids obtained by ``dragging'' a chosen puncture $x$ of $R$ on $S$ about an
embedded loop on $S$ that avoids the remaining punctures of $R$ on $S$. These particular bounding pair braids generate a normal subgroup
of $PMod_R(S)$ which we identify naturally with the fundamental group $\pi_1(R \cup \{x\}, x)$ of $R \cup \{x\}$ based at $x$. We show
that an arbitrary automorphism of $PMod_R(S)$ maps $\pi_1(R \cup \{x\}, x)$ isomorphically to $\pi_1(R \cup \{y\}$ for some puncture $y$
of $R$ on $S$. By pulling back a given automorphism by an automorphism induced by a homeomorphism $R \rightarrow R$, we are able to
assume that $y = x$. In this way, we reduce the problem to studying automorphisms of $PMod_R(S)$ which restrict to automorphisms of the
surface group $\pi_1(R \cup \{x\}, x)$. 

It is a well-known fact that an automorphism of $\pi_1(R \cup \{x\}, x)$ is induced by a homeomorphism $(R \cup \{x\}, x) \rightarrow (R
\cup \{x\}, x)$ if and only if it preserves the peripheral structure of $\pi_1(R \cup \{x\}, x)$. We show that this peripheral structure
cooresponds to the $2$-string twist braids on $R$ which have one of their two strings based at $x$. From this we deduce that
automorphisms of $\pi_1(R \cup \{x\}, x)$ which are induced, via pull-back by homeomorphisms and restriction, from automorphisms of
$PMod_R(S)$ preserve the peripheral structure of $\pi_1(R \cup \{x\}, x)$ and, hence, are induced by homeomorphisms 
$(R \cup \{x\}, x) \rightarrow (R \cup \{x\}, x)$. 

This allows us to reduce the problem to studying automorphisms of $PMod_R(S)$
which are the identity on the surface subgroup $\pi_1(R \cup \{x\}, x)$. We show
that any such automorphism of $PMod_R(S)$ is the identity automorphism. This
proves Theorem $2$. 

Here is an outline of the paper. In Section \ref{sec:preliminaries}, we review the
basic notions and results related to mapping class groups and surface braid
groups. We assume that the reader is familiar with the fundamentals of Thurston's
theory of surfaces (cf. \cite{flp}) 

In Section \ref{sec:pureispure}, we prove that elements of $PMod_R(S)$ are pure in the sense defined by Ivanov \cite{i2}. In Section
\ref{sec:ranks}, we use the results of Section \ref{sec:pureispure} to develop an explicit formula for the rank of the center of the centralizer of
an element of $PMod_R(S)$ in terms of its Thurston normal forms. 

In Section \ref{sec:basicbraids}, we use the rank formula of Section \ref{sec:ranks} to determine which braids are basic braids, pure
braids which have rank one centers of centralizers. In Section \ref{sec:bpbraids}, we give an algebraic characterization of bounding pair braids,
thereby distinguishing them algebraically from the other two types of basic braids, twist braids and basic single-$pA$ braids. In
Section \ref{sec:twistbraids}, we algebraically distinguish twist braids from basic single-$pA$ braids. In Section \ref{sec:twiststrings}, we
algebraically characterize the number of strings in a twist braid. In Section \ref{sec:bpstrings}, we algebraically characterize the number of
strings in a bounding pair braid. 

In Section \ref{sec:punctures}, we algebraically characterize the situation when two $1$-string bounding pair braids are based at the same puncture
of $R$ on $S$. In Section \ref{sec:surfsubgp}, we use the results of Section \ref{sec:punctures} to prove that automorphisms of $PMod_R(S)$
preserve the surface subgroups $\pi_1(R \cup \{x\}, x)$. In Section \ref{sec:peripheral}, we prove that the resulting induced automorphisms of the
surface subgroups $\pi_1(R \cup \{x\}, x)$ preserve the peripheral structure and are, hence, {\it geometric} (i.e induced by a homeomorphism $R
\rightarrow R$). 

In Section \ref{sec:mainresults}, we assemble the results of the previous sections to prove the main results of this paper, Theorems $1$ and $2$.

\section{Preliminaries} 
\label{sec:preliminaries} 

Let $C(R)$ denote the complex of curves of $R$. Let $\alpha$ be a vertex of $C(R)$. We say that $\alpha$ is
$S$-essential if each representative circle $C$ of $\alpha$ is essential on $S$, as well as on $R$.
Otherwise, we say that $\alpha$ is $S$-inessential. Note that $\alpha$ is $S$-inessential if and only if
each representative circle $C$ of $\alpha$ bounds a disc $D$ with $k$ punctures where $2 \leq k \leq n$.
Note, furthermore, that the integer $k$ does not depend upon the representative circle $C$ of $\alpha$. 

Let $f$ be a twist about an essential circle $C$ on $R$. Note that $f$ is a nontrivial element of $PMod_R$. Moreover, 
$f$ is an element of the pure braid group, $PMod_R(S)$, if and only if $C$ is not essential on $S$. If $C$ bounds a
$k$-punctured disc $D$ on $S$ (i.e. an embedded disc $D$ on $S$ with its boundary contained in $R$ and with exactly $k$ 
punctures of $R$ in its interior) and $p_i$, $1 \leq i \leq k$, are the $k$ punctures of $R$ in the disc $D$, then we say that $f$ is a
{\it $k$-string twist braid supported on the $k$-punctured disc $D$ with strings based at $p_i$, $1 \leq i \leq k$}.

Suppose that $C$ and $D$ are disjoint, non-isotopic, essential circles on $R$ which are also essential on $S$. Let $f$ be the
product $t_C \circ t_D^{-1}$ of opposite twists about $C$ and $D$. Note that $f$ is a nontrivial element of $PMod_R$. Moreover,
$f$ is an element of the pure braid group, $PMod_R(S)$, if and only if $C$ and $D$ are isotopic on  $S$. If $C$ and $D$ co-bound
a $k$-punctured annulus $A$ on $S$ (i.e. an embedded annulus $A$ on $S$ with its boundary contained in $R$ and with exactly $k$ 
punctures of $R$ in its interior) and $p_i$, $1 \leq i \leq k$, are the $k$ punctures of $R$ in the annulus $A$, , then we say that $f$
is a {\it $k$-string bounding pair braid supported on the $k$-punctured annulus $A$ with strings based at $p_i$, $1 \leq i \leq k$}. 

Suppose that $k = 1$. Let $x$ be the unique puncture of $R$ which is contained in the annulus $A$. Equip $A$ and its boundary
$C \cup D$ with the orientation induced from the orientation of the surface $S$. Up to isotopy, there exists a unique loop $c$ in
the interior of $A$ based at $x$ such that $c$ is isotopic to the oriented circle $C$ on $S$. Note that, if $n = 1$, the
$1$-string bounding pair braid $t_C \circ t_D^{-1}$ supported on the annulus $A$ is the spin of $x$ about $c$ in the sense of
\cite{b}. We shall say that $t_C \circ t_D^{-1}$ is the $1$-string bounding pair braid on $R$ corresponding to the embedded 
loop $c$ in $R \cup \{x\}$ based at $x$. 

\section{Pure braids are pure} 
\label{sec:pureispure}

In this section, we prove that pure surface braids are pure in the sense defined by Ivanov, \cite{i2}. 

\begin{thm} Let $f$ be an element of the pure braid group $PMod_R(S)$. Let $C$ be a system of circles on $R$ representing the essential reduction system
$\mathcal{C}$ of $f$. Let $C'$ be the union of all components of $C$ which are essential on $S$. There exists a homeomorphism $F : R \rightarrow R$ such
that: 
\begin{itemize} 
\item[(i)] $F$ fixes pointwise each element of $C$, 
\item[(ii)] $F$ maps each component of the complement of $C$ on $R$ to itself, 
\item[(iii)] the restriction of $F$ to each component of the complement of $C$ on $R$ is either pseudo-Anosov or isotopic to the identity. 
\end{itemize} 

\noindent and if $\overline{F} : S \rightarrow S$ is the extension of $F$ to $S$, then:

\begin{itemize} 
\item[(iv)] there exists an isotopy on $S$, preserving each component of $C'$, from the identity map of $S$ to $\overline{F}$.
\end{itemize}

\label{thm:pureispure} \end{thm} 

\section{ranks of centers of centralizers of pure braids} 
\label{sec:ranks}

Let $f$ be an element of the pure braid group $PMod_R(S)$. 

Let $\sigma \subset C(R)$ be the essential reduction system of $f$. 

Let $\sigma^*$ be the set of $S$-essential elements of $\sigma$ (i.e. $\sigma^* = \sigma \cap C^*(R)$). For each element
$\alpha$ of $\sigma^*$, let $[\alpha]$ be the set of all elements of $\sigma^*$ which are $S$-equivalent to $\alpha$.
Then there exist a nonnegative integer $k$ and a subset ${\alpha_i | 1 \leq j \leq k}$ of $\sigma^*$ such
that $\sigma*$ is the disjoint union of $[\alpha_j]$, $1 \leq j \leq k$. Note that $\sigma^*$ is empty if
and only if $k = 0$. 

Let $p$ be the number of $pA$ components of the reduction of $f$ along $\sigma$. Let $i$ be the
number of $S$-inessential elements of $\sigma$. For each integer $j$ such that $1 \leq j \leq k$, let $e_j$
denote the number of elements of $[\alpha_j]$. 

\begin{prop} The rank $\mu$ of the center of the centralizer of $f$ is equal to $p + i + \Sigma_{1 \leq j
\leq k} (e_j - 1)$. 
\label{prop:rank}
\end{prop}

\begin{proof} This result is an application of Theorem \ref{thm:pureispure}. 
\end{proof}

\section{basic braids}
\label{sec:basicbraids} 

Let $f$ be an element of the pure braid group $PMod_R(S)$. We say that $f$ is a {\it basic braid} 
if the rank $\mu$ of the center of the centralizer of $f$ in $PMod_R(S)$ is equal to $1$. In this section, we shall give a
detailed description of basic braids. 

Let $f$ be an element of the mapping class group $Mod_R$ of $R$. 

We say that $f$ is a {\it single-$pA$} mapping class if (i) the reduction of $f$ along its essential reduction system $\sigma$
has exactly one $pA$ component $C$ and (ii) the essential reduction system $\sigma$ of $f$ corresponds to the boundary of this
component $C$ (i.e. $\sigma$ is equal to the set of isotopy classes of boundary components of $C$ which are essential on $R$).  

\begin{prop}
\label{prop:basic} Let $f$ be an element of the braid group $PMod_R(S)$ such that the center of the centralizer
of $f$ in $PMod_R(S)$ is infinite cyclic. Then one of the following holds:

(i) $f$ is a single-$pA$ mapping class supported on a subsurface $T$ of $R$ such that each boundary component
of $T$ is essential on $S$ and no two boundary components of $T$ are isotopic on $S$.

(ii) $f$ is a $k$-string bounding pair braid for some integer $k$ with $1 \leq k \leq n$.

(iii) $f$ is a $k$-string twist braid for some integer $k$ with $2 \leq k \leq n$.
\end{prop}

Hence, the basic braids are precisely the braids described in $(i)$, $(ii)$ and $(iii)$ of Proposition \ref{prop:basic}. We
shall refer to the braids described in $(i)$ as basic single-$pA$ braids. We shall refer to the braids described in $(i)$ as
basic $S$-bounding pair braids. Note that the braids described in $(ii)$ are, according to the terminology of \cite{b}, {\it
spins}. We shall refer to the braids described in $(iii)$ as basic {\it twist} braids.   

\begin{rem} If $f$ is a power of the product of opposite twists on two disjoint, non-isotopic circles in $R$
which co-bound a once-punctured annulus in $R$, and one of these circles is inessential in $S$, then both of
these circles will be inessential in $S$, and the center of the centralizer of $f$ in $PMod_R(S)$ will be the
free abelian subgroup of rank $2$ in $PMod_R(S)$ freely generated by the twists about these two circles.
\end{rem}

\section{bounding pair braids}
\label{sec:bpbraids} 

In this section, we shall distinguish basic bounding pair braids algebraically from basic single-$pA$ braids and basic twist braids when $R$
has at least $2$ punctures.

\begin{prop} Let $n \geq 2$. Let $f$ be an element of the pure braid group $PMod_R(S)$ such that the center of the centralizer
of $f$ in $PMod_R(S)$ is infinite cyclic. Then the following are equivalent: 

\noindent (a) $f$ is a bounding pair braid

\noindent (b) There exists a basic braid $g$ in $PMod_R(S)$ such that $f$ and $g$ generate a free abelian group of rank
$2$ in $PMod_R(S)$ and $fg$ is a basic braid.
\label{prop:bpbraids}
\end{prop}

\begin{rem} Proposition \ref{prop:bpbraids} holds for all closed surfaces $S$ and any number of punctures $n$. 
\end{rem}

\begin{rem} Proposition \ref{prop:bpbraids} is vacuously true when the genus of $S$ is zero. After all, there are no
bounding pair braids on surfaces of genus zero. Indeed, there are no essential simple closed curves on a closed surface of
genus zero. Moreover, if the genus of $S$ is zero, then there are no pairs of basic braids as in clause $(b)$ of Proposition 
\ref{prop:bpbraids}.\end{rem}  

\section{twist braids} 
\label{sec:twistbraids} 

In this section, we shall distinguish twist braids algebraically from basic single-$pA$ braids when $R$ has at least $2$ punctures.
  
This will require the following algebraic concept. Let $A$ and $B$ be elements of a group $G$. Let $H$ be the subgroup of $G$
generated by $A$ and $B$. We say that $A$ adheres to $B$ in $G$ if $A$ commutes with $D$, whenever $D$ is an element of the
centralizer of some nontrivial element $C$ of $H$.

Using this algebraic concept of adherence, we can distinguish twist braids algebraically from basic single-$pA$ braids as follows:

\begin{prop} Let $f$ be an element of the pure braid group $PMod_R(S)$. Suppose that the center of the
centralizer of $f$ in $PMod_R(S)$ is infinite cyclic. Suppose that $f$ is not a bounding pair braid. Then the following
are equivalent:

\noindent (a) $f$ is a twist braid (i.e. $f$ is not a single-$pA$ braid).

\noindent (b) There exists an element $g$ in $PMod_R(S)$ such that $f$ and $g$ generate a free abelian group of rank $2$ in
$PMod_R(S)$ and $f$ adheres to $g$ in $PMod_R(S)$.
\label{prop:twistbraids} \end{prop}

\begin{rem} Since the property expressed in part (b) of Proposition \ref{prop:twistbraids} is an algebraic property, we can use
this property to distinguish twist braids algebraically from basic single-$pA$ braids.
\end{rem}

\section{counting strings in twist braids}
\label{sec:twiststrings} 

In the previous sections, we algebraically distinguished the three types of basic braids occurring in Proposition \ref{prop:basic},
basic single-$pA$ braids, bounding pair braids, and twist braids. In this section, we shall show how to algebraically detect the number
of strings in a twist braid. 

First, we show how to algebraically detect $2$-string twist braids. 

\begin{prop} Let $f$ be an element of the pure braid group $PMod_R(S)$. Suppose that $f$ is a $k$-string twist braid for
some integer $k$ with $2 \leq k \leq n$. Then the following are equivalent: 

\noindent (a) $k = 2$

\noindent (b) There exist basic bounding pair braids $g_j$, $1 \leq j \leq n - 2$ such that $f$ and $g_j$, $1 \leq j \leq n -
2$ generate a free abelian group of rank $n - 1$ in $PMod_R(S)$. 
\label{prop:2twist} \end{prop} 

Now, we may algebraically detect $k$-string twist braids inductively as follows. 

\begin{prop} Let $f$ be an element of the pure braid group $PMod_R(S)$. Suppose that $k$ is an integer such that $3 \leq k \leq
n$. Suppose that $f$ is a twist braid. Then the following are equivalent: 

\noindent (a) $f$ is a $k$-string twist braid. 

\noindent (b) $f$ is not an $l$-string twist braid for any integer $l$ less than $k$ and there exists bounding pair braids, 
$g_j$, $1 \leq j \leq n - k$ such that $f$ and $g_j$, $1 \leq j \leq n - k$ generate a free abelian group of rank $n + 1 - k$
in $PMod_R(S)$.
\label{prop:ktwist} \end{prop}

\section{counting strings in bounding pair braids}
\label{sec:bpstrings} 

In the previous section, we showed how to algebraically detect the number of strings in a twist braid. In this section, we shall
use this to algebraically detect the number of strings in a bounding pair braid.

Throughout this section, we shall use the following notation. Let $f$ be a $k$-string bounding pair braid. Let $C$ and $D$ be
a pair of circles on $R$ associated to $f$. By definition, $f$ is the product of opposite twists about $C$ and $D$ and $C$ and $D$
bound an annulus $A$ on $S$ such that $A$ contains exactly $k$ punctures of $R$. 

Note that the complement of $A$ in $S$ has one or two components, depending upon whether $C$ and $D$ are both nonseparating or
both separating circles on $S$. By the definition of a bounding pair braid, neither $C$ nor $D$ bounds a disc. Hence, each
component of the complement of $A$ in $S$ has positive genus.

We say that a component $K$ of the complement of $C \cup D$ in $S$ is {\it adequately punctured} if it contains at least two
punctures of $R$. We begin with showing how to algebraically detect adequately punctured components of the complement of 
$C \cup D$ in $S$. To do this, we consider the collection $\Delta (f)$ of twist braids $\alpha$ such that $\alpha$ commutes with $f$.
Note that this collection corresponds to the collection of isotopy classes of essential circles $A$ on $R$ which are disjoint from
both $C$ and $D$ and bound discs on $S$ which contain at least two punctures of $R$. 

\begin{rem} Note that $\Delta (f)$ is nonempty unless either (i) $n = 2$ and the complement of $C \cup D$ in $S$ has exactly two
components each with exactly one puncture or (ii) $n = 3$ and the complement of $C \cup D$ in $S$ has exactly three components
each with exactly one puncture. In all other situations, at least one of the components of the comlement of $C \cup D$ in $S$ is
adequately punctured and, hence, supports a twist braid $\alpha$ commuting with $f$. 
\end{rem} 

\begin{prop} Let $f$ be a $k$-string bounding pair braid. Let $C$ and $D$ be a pair of circles on $R$ associated to $f$. Let
$\alpha$ and $\beta$ be elements of $\Delta (f)$. Let $A$ and $B$ be circles on $R$ associated to $\alpha$ and $\beta$ such that $A$
and $B$ are both disjoint from $C$ and $D$. Then the following are equivalent: 

\noindent (a) There exists an element $\gamma$ in $\Delta (f)$ such that $\alpha$ and $\beta$ do not commute with $\gamma$.

\noindent (b) $A$ and $B$ are contained in the same component of the complement of $C \cup D$ in $S$. 
\label{prop:equivreln} \end{prop} 

It follows that the relation $\sim$ on $\Delta (f)$ defined by part (a) of Proposition \ref{prop:equivreln} is an equivalence
relation on $\Delta (f)$. 

\begin{prop} Let $\sim$ be the relation on $\Delta (f)$ defined by the rule $\alpha \sim \beta$ if and only if there exists an
element $\gamma$ in $\Delta (f)$ such that $\alpha$ and $\beta$ do not commute with $\gamma$. Then $\sim$ is an equivalence
relation on $\Delta (f)$. 
\label{prop:simequivreln} \end{prop} 

Moreover, by Proposition \ref{prop:equivreln}, it follows that the adequately punctured components of the complement of $C
\cup D$ in $S$ are in natural one-to-one correspondence with the equivalence classes of the equivalence relation $\sim$ on
$\Delta (f)$. 

\begin{prop} There exists a well-defined bijection $\beta : \Delta (f)/\sim \rightarrow \pi_0(S \setminus A)$ such that
$\beta(\alpha)$ is equal to the component $K$ of $S \setminus A$ if and only if $K$ contains a circle $A$ which is associated to
$\alpha$ and is disjoint from $C$ and $D$.
\label{prop:bijectionbeta} \end{prop}

\begin{cor} The number of adequately punctured components of the complement of $C \cup D$ in $S$ is equal to the number of
equivalence classes of the equivalence relation $\sim$ on $\Delta (f)$.
\label{cor:adeqpunctured} \end{cor} 

Next, we show how to algebraically detect the number of punctures in an adequately punctured component $K$ of the complement of
$C \cup D$ in $S$.

\begin{prop} Let $K$ be an adequately punctured component of the complement of $C \cup D$ in $S$. The number of punctures of
$R$ contained in $K$ is equal to the maximum number of strings of a twist braid in $\chi^{-1}(K)$.  
\label{prop:countingpunctures} 
\end{prop}

Next, we show how to algebraically detect an adequately punctured component of positive genus. 

\begin{prop} Let $K$ be a $k$-punctured component of the complement of $C \cup D$ in $S$ where $2 \leq k \leq n$. Then the
following are equivalent: 
 
\noindent (a) The genus of $K$ is positive.  

\noindent (b) For each $k$-twist $\alpha$ supported in $K$ and each $(k - 1)$-twist $\beta$ supported in $A$ such that $\alpha$
and $\beta$ commute, there exists an element $g$ in the pure braid group $PMod_R(S)$ such that (i) $g$ commutes with $f$ and 
$\beta$ and (ii) $f$ adheres to $g$.

\label{prop:positivegenus} 
\end{prop}

Note that the above results, taken together, algebraically characterize the number of strings in a bounding pair braid. 

\section{detecting punctures}
\label{sec:punctures} 

In this section, we show how to algebraically detect punctures. 

Let $x$ be a puncture of $R$. Let $\mathcal{N}(x)$ denote the collection of all $1$-string bounding pair braids $f$
corresponding to embedded loops $c$ in $R \cup \{x\}$ based at $x$. 

\begin{prop} Let $x$ and $y$ be punctures of $R$. Let $f$ be an element of $\mathcal{N}(x)$. Let $g$ be an element of
$\mathcal{N}(y)$. Then the following are equivalent: 

\noindent (a) $x = y$

\noindent (b) There exists a sequence $f_i$, $1 \leq i \leq N$, of $1$-string bounding pair braids $f_i$ such that 
(i) $f_1 = f$, (ii) $f_N = g$, and (iii) for each integer $i$ with $1 \leq i < N$ there exists an $(n - 1)$-string twist braid
$t_i$ such that $f_i$ and $f_{i + 1}$ both commute with $t_i$. 
\label{prop:detectpunctures} \end{prop} 

\section{preservation of surface subgroups} 
\label{sec:surfsubgp}

In this section, we shall show that any automorphism of the pure braid group $PMod_R(S)$ induces an isomorphism between certain
natural normal subgroups of $PMod_R(S)$ which are naturally isomorphic to fundamental groups of surfaces.

Let $x$ be a puncture of $R$ on $S$. Note that there is a natural homomorphism $\mu: PMod_R \rightarrow PMod_{R \cup \{x\}}$
corresponding to extending homeomorphisms $F : R \rightarrow R$ which ``fix'' each puncture of $R$ on $S$ to homeomorphisms 
$\overline{F} : R \cup \{x\} \rightarrow R \cup \{x\}$. Following the notation $PMod_R(S)$ introduced above, we shall denote the
kernel of $\mu$ as $PMod_R(R \cup \{x\})$. Note that $PMod_R(R \cup \{x\})$ is a subgroup of $PMod_R(S)$. Since $PMod_R(R \cup \{x\})$ is
a normal subgroup of $PMod_R$, it is also a normal subgroup of the subgroup $PMod_R(S)$ of $PMod_R$ in which it is contained.  

Recall that there is a natural homomorphism $\eta: \pi_1(R \cup \{x\}, x) \rightarrow PMod_R(S)$ which sends the element $[c]$ of
$\pi_1(R \cup \{x\}$ corresponding to any embedded loop $c$ in $R \cup \{x\}$ based at
$x$ to the corresponding $1$-string bounding pair braid $t_C \circ t_D^{-1}$ on $R$, where $C \cup D$ is the boundary of a closed annular
neighborhood $A$ of $c$ on $S$ such that $C \cup D$ is contained in $R$ and $x$ is the unique puncture of $R$ in the interior of $A$.  
This homomorphism arises from the long exact homotopy sequence of the evaluation map $ev_x : Diff(R \cup \{x\}) \rightarrow R \cup \{x\}$
defined by the rule $ev_x(F) = F(x)$. $\eta$ is injective provided $R \cup \{x\}$ has positive genus (cf. Theorems 1.4, 4.2, and 4.3 of
\cite{b}).  

\begin{prop} Let $x$ be a puncture of $R$. The natural monomorphism 
$\eta: \pi_1(R \cup \{x\}, x) \rightarrow PMod_R(S)$ maps $\pi_1(R \cup \{x\}, x)$ isomorphically onto $PMod_R(R \cup \{x\})$.
\label{prop:naturalisomorphism} \end{prop} 

Using the natural monomorphism $\eta: \pi_1(R \cup \{x\}, x) \rightarrow PMod_R(S)$, we shall identify 
$\pi_1(R \cup \{x\}, x)$ with its image under $\eta$. This identification identifies the natural action of homeomorphisms $\overline{G} :
(R \cup \{x\}, x) \rightarrow (R \cup \{x\}, x)$ of the pointed space $(R \cup \{x\}, x)$ on $\pi_1(R \cup \{x\}, x)$ with the
natural action of their restrictions $G : R \rightarrow R$ on $PMod_R(R \cup \{x\})$. 

\begin{prop} Let $\overline{G} : (R \cup \{x\}, x) \rightarrow (R \cup \{x\}, x)$ be a homeomorphism of the pointed space 
$(R \cup \{x\}, x)$. Let $G : R \rightarrow R$ be the restriction of $\overline{G}$ to $R$. Let $G_*: PMod_R(S) \rightarrow PMod_R(S)$
be the automorphism of $PMod_R(S)$ defined by the rule $G_*([H]) = [G \circ H \circ G^{-1}]$ for each mapping class $[H]$ in
$PMod_R(S)$. Then $G_*(\eta(a)) = \eta(\overline{G}_*(a))$ for each element $a$ in $\pi_1(R \cup \{x\}, x)$. 
\label{prop:equivariance} \end{prop} 

The following proposition demonstrates the significance of the surface subgroups $\pi_1(R \cup \{x\}, x)$ for our purposes. 

\begin{prop} Let $x$ be a puncture of $R$ on $S$. Let $\phi : PMod_R(S) \rightarrow PMod_R(S)$ be an automorphism of $PMod_R(S)$.
If $\phi$ fixes each element of the subgroup $\pi_1(R \cup \{x\}, x)$ of $PMod_R(S)$, then $\phi$ is the identity automorphism of
$PMod_R(S)$. 
\label{prop:surfsubgp} \end{prop} 

\begin{proof} Let $f$ be an element of $PMod_R(S)$. Let $h = f^{-1} \phi(f)$. Note that $h$ is an element of $PMod_R(S)$. Hence, $h$ is 
equal to the mapping class $[H]$ of some orientation preserving homeomorphism $H : R \rightarrow R$ of $R$. We shall show that $H$ is
isotopic to the identity. 

Let $C$ be an essential circle on $R$. Choose a closed annulus $A$ on $S$ such that the boundary of $A$ is contained in
$R$, $C$ is one of the two boundary components of $A$, and $x$ is the unique puncture of $R$ on $S$ contained in the interior of $A$.
Let $c$ be an embedded loop in the interior of $A$ based at $x$ such that $c$ is isotopic to $C$. By the previous observations, the
element $[c]$ of $\pi_1(R \cup \{x\}, x)$ represented by $c$ corresponds to the $1$-string bounding pair braid $t_C \circ t_D^{-1}$ on
$R$ where $D$ is the boundary component of $A$ distinct from $C$.  

Let $g$ denote $t_C \circ t_D^{-1}$. Since $g$ is an element of $\pi_1(R \cup \{x\}, x)$, $\phi(g) = g$. 

Let $f$ be an element of $PMod_R(S)$. Since $\pi_1(R \cup \{x\}, x)$ is a normal subgroup of $PMod_R(S)$ and $g$ is an element of
$\pi_1(R \cup \{x\}, x)$, it follows that $fgf^{-1}$ is also an element of $\pi_1(R \cup \{x\}, x)$. Hence, $\phi(fgf^{-1}) =
fgf^{-1}$. On the other hand, since $\phi$ is a homomorphism on $PMod_R(S)$, $\phi(fgf^{-1}) = \phi(f) \phi(g) \phi(f)^{-1}$. Hence,
$fgf^{-1} = \phi(f) \phi(g) \phi(f)^{-1} = \phi(f) g \phi(f)^{-1}$. 

Since $g = t_C \circ t_D^{-1}$ and $H : R \rightarrow R$, it follows that $hgh^{-1} = [H] \circ t_C \circ t_D^{-1} \circ [H]^{-1} =
t_{H(C)} \circ t_{H(D)}^{-1}$. Since $hgh^{-1} = g$, this implies that $t_{H(C)} \circ t_{H(D)}^{-1} = t_C \circ t_D^{-1}$. It follows,
by standard arguments, that $H(C)$ is isotopic to $C$, and $H(D)$ is isotopic to $D$.

This proves that $H$ preserves the isotopy class of every essential circle $C$ on $S$. It follows by Lemma 5.1 and Theorem 5.3 of
\cite{im}, that $H$ is isotopic to the identity. In other words, the mapping class $h$ of $H$ is the trivial element $id$ of $PMod_R(S)$.
Thus, $f^{-1} \phi(f) = id$ and, hence, $\phi(f) = f$. 

We have shown that $\phi(f) = f$ for each element $f$ of $PMod_R(S)$. Hence, $\phi$ is the identity automorphism of $PMod_R(S)$. 
\end{proof}

Motivated by Proposition \ref{prop:surfsubgp}, we now turn to an investigation of the restriction of an arbitrary
automorphism of $PMod_R(S)$ to a surface subgroup $\pi_1(R \cup \{x\}, x)$ of $PMod_R(S)$.

\begin{prop} Let $\chi : PMod_R(S) \rightarrow PMod_R(S)$ be an automorphism of the pure braid group $PMod_R(S)$. Let $x$ be a puncture
of $R$ on $S$.  There exists a unique puncture $y$ of $R$ on $S$ such that $\chi$ maps $\pi_1(R \cup \{x\}, x)$ to $\pi_1(R \cup
\{y\},y)$.
\label{prop:pressurfsubgroups} \end{prop}
 
\begin{prop} Let $\chi : PMod_R(S) \rightarrow PMod_R(S)$ be an automorphism of the pure braid group $PMod_R(S)$. Let $x$ be a puncture
of $R$ on $S$.  Let $y$ be the unique puncture of $R$ on $S$ such that $\chi$ maps $\pi_1(R \cup \{x\}, x)$ to $\pi_1(R \cup
\{y\},y)$. Then the restriction $\chi|: \pi_1(R \cup \{x\}, x) \rightarrow \pi_1(R \cup \{y\}, y)$ is an isomorphism.
\label{prop:inducedisomorphism} \end{prop}
 
\begin{proof} Let $\tau : PMod_R(S) \rightarrow PMod_R(S)$ be the inverse of $\chi : PMod_R(S) \rightarrow PMod_R(S)$. Let $z$ be the
unique puncture of $R$ on $S$ such that $\tau$ maps $\pi_1(R \cup \{y\}, y)$ to $\pi_1(R \cup \{z\},z)$. Note that the composition 
$\tau \circ \chi$ maps $\pi_1(R \cup \{x\}, x)$ into $\pi_1(R \cup \{z\}, z)$. On the other hand, this composition is the identity map
of $PMod_R(S)$. Hence, $\pi_1(R \cup \{x\}, x)$ is a subgroup of $\pi_1(R \cup \{z\}, z)$ (i.e. the image of $\pi_1(R \cup \{x\}, x)$
under the natural monomorphism $\eta_x : \pi_1(R \cup \{x\}, x) \rightarrow PMod_R(S)$ is a subgroup of the image of $\pi_1(R \cup
\{z\}, z)$ under the natural monomorphism $\eta_z : \pi_1(R \cup \{z\}, z) \rightarrow PMod_R(S)$. 

In particular, every $1$-string bounding pair braid $f$ with strings based at $x$ is an element of $\pi_1(R \cup \{z\}, z)$. Since
$\pi_1(R \cup \{z\}, z) \subset PMod_R(S)$ is the kernel of the natural homomorphism 
$ext_z: PMod_R(S) \rightarrow PMod_{R \cup \{z\}}(S)$, this implies that $ext_z(f)$ is the trivial element of $PMod_{R \cup \{z\}}(S)$. 

Suppose that $z$ is not equal to $x$. Then, clearly, $ext_z(f)$ is a $1$-string bounding pair braid on $R \cup \{z\}$ with string based
at $x$ on $R \cup \{z\}$. Hence, $ext_z(f)$ is not the trivial element of $PMod_{R \cup \{z\}}(S)$. This is a contradiction. Hence, 
$z = x$. 

Since $z = x$, $\tau$ maps $\pi_1(R \cup \{y\}, y)$ to $\pi_1(R \cup \{x\},x)$. Clearly, since $\chi$ and $\tau$ are inverse
homomorphisms, the restrictions $\chi| : \pi_1(R \cup \{x\}, x) \rightarrow \pi_1(R \cup \{y\}, y)$ and $\tau| : \pi_1(R \cup \{y\}, y)
\rightarrow \pi_1(R \cup \{x\}, x)$ are inverse homomorphisms. In particular, $\chi| : \pi_1(R \cup \{x\}, x) \rightarrow \pi_1(R \cup
\{y\}, y)$ is an isomorphism. 
\end{proof} 

\begin{prop} Let $\chi : PMod_R(S) \rightarrow PMod_R(S)$ be an automorphism of the pure braid group $PMod_R(S)$. Let $x$ be a puncture
of $R$ on $S$.  Let $y$ be the unique puncture of $R$ on $S$ such that $\chi$ maps $\pi_1(R \cup \{x\}, x)$ to $\pi_1(R \cup
\{y\},y)$. Let $F : R \rightarrow R$ be the restriction of a homeomorphism $\overline{F} : (R \cup \{y\}, y) \rightarrow (R \cup \{x\},
x)$ of pointed spaces. Let $F_* : PMod_R(S) \rightarrow PMod_R(S)$ be the automorphism of $PMod_R(S)$ defined by the rule $F_*([H]) = [F
\circ H \circ F^{-1}]$ for evey mapping class $[H]$ in $PMod_R(S)$. Let $\theta : PMod_R(S) \rightarrow PMod_R(S)$ be the composition 
$F_* \circ \chi : PMod_R(S) \rightarrow PMod_R(S)$. Then $\theta$ is an automorphism of $PMod_R(S)$ and $\theta$
restricts to an automorphism $\theta| : \pi_1(R \cup \{x\}, x) \rightarrow \pi_1(R \cup \{x\}, x)$ of $\pi_1(R \cup \{x\}, x)$. 
\label{prop:inducedautomorphism} \end{prop} 

\begin{proof} Clearly, $F_* : PMod_R(S) \rightarrow PMod_R(S)$ restricts to an isomorphism, $F_*|: \pi_1(R \cup \{x\}, x) \rightarrow
\pi_1(R \cup \{y\}, y)$. Since $\chi : PMod_R(S) \rightarrow PMod_R(S)$ also restricts to an isomorphism, $\chi|: \pi_1(R \cup \{x\}, x)
\rightarrow \pi_1(R \cup \{y\}, y)$, the composition $\theta: PMod_R(S) \rightarrow PMod_R(S)$ of $F_*$ and $\chi$ restricts to
an automorphism, $\theta|: \pi_1(R \cup \{x\}, x) \rightarrow \pi_1(R \cup \{x\}, x)$, of $\pi_1(R \cup \{x\}, x)$. 
\end{proof}

\section{preservation of peripheral structures}
\label{sec:peripheral}

In this section, we shall prove that the restrictions of automorphisms of the pure braid group, $PMod_R(S)$, to surface subgroups, as in
Proposition \ref{prop:inducedisomorphism}, respect the peripheral structure of these subgroups.

Let $\chi$ be an automorphism of $PMod_R(S)$. As in Proposition \ref{prop:inducedautomorphism}, we
may choose a homeomorphism $F : R \rightarrow R$ of $R$ such that the composition $\theta$ of $F_*$ and $\chi$ restricts to an
automorphism $\theta|$ of $\pi_1(R \cup \{x\}, x)$. Speaking more precisely about what we shall do in this section, we shall prove that
the restriction $\theta|: \pi_1(R \cup \{x\}, x) \rightarrow \pi_1(R \cup \{x\}, x)$ respects the peripheral structure of $\pi_1(R \cup
\{x\}, x)$. This is the precise condition for ensuring that this restriction of $\theta$ is induced by a homeomorphism of $R$ sending
$x$ to $x$. As we shall see later, this is enough to ensure that $\theta$ is induced by a homeomorphism of $R$. From this, it will easily
follow that $\chi$ is induced by a homeomorphism of $R$.  

Due to the naturality of the monomorphism $\eta: \pi_1(R \cup \{x\}, x) \rightarrow PMod_R(S)$, the desired result
concerning the restriction $\theta|$ can be achieved by appealing to a classical result about automorphisms of fundamental groups of
punctured surfaces. Namely, an automorphism $\phi$ of $\pi_1(R \cup \{x\}, x)$ is induced by ahomeomorphism $G : (R \cup
\{x\}, x) \rightarrow (R \cup \{x\}, x)$ if and only if $\phi$ preserves the peripheral structure of $\pi_(R \cup \{x\}, x)$. 
Our task, therefore, will be to show that the restriction $\theta|$ preserves the peripheral structure of $\pi_1(R \cup \{x\}, x)$.
Hence, we need to understand the peripheral structure of $\pi_1(R \cup \{x\}, x)$ and its relationship to the natural identification of
$\pi_1(R \cup \{x\}, x)$ with a subgroup of $PMod_R(S)$ via the natural monomorphism $\eta$. 

Let $x_i$, $1 \leq i \leq n$, be the $n$ punctures of $R$. Choose $n$ disjoint embedded discs, $D_i$, $1 \leq i \leq n$, on $S$
such that $x_i$ lies in the interior $U_i$ of $D_i$ in $S$. Let $T_i$ denote the compact subsurface of $S$ obtained by deleting
the interiors $U_j$ of each of the discs $D_j$ with $1 \leq j \leq n$ and $j \neq i$. Note that $x_i \in T_i \subset R \cup
\{x_i\}$ and $(T_i, x_i)$ is a deformation retract of $(R \cup \{x_i\}, x_i)$. In particular, the inclusion $\imath : (T_i, x_i)
\rightarrow (R \cup \{x_i\}, x_i)$ induces an isomorphism $\imath_*: \pi_1(T_i, x_i) \rightarrow \pi_1(R \cup \{x_i\}, x_i)$ of
fundamental groups. Recall that, by definition, the peripheral structure of $\pi(T_i, x_i)$ is the set $P_i$ of $2(n - 1)$ conjugacy
classes in $\pi_1(T_i, x_i)$ corresponding to the oriented boundary components of $T_i$. We shall refer to the corresponding subset 
$\imath_*(P_i)$ of $\pi_1(R \cup \{x_i\}, x_i)$ as the peripheral structure of $\pi_1(R \cup \{x_i\}, x_i)$ and denote this subset as 
$\mathcal{P}_i$. 

\begin{prop} The image of the peripheral structure $\mathcal{P}_i$ of $\pi_1(R \cup \{x_i\}, x_i)$ under the natural monomorphism 
$\eta: \pi_1(R \cup \{x_i\}, x_i) \rightarrow PMod_R(S)$ is equal to the set of all $2$-string twist braids $f$ such that there
exists an integer $j$ with $1 \leq j \leq n$ and $j \neq i$ such that $f$ has strings based at $x_i$ and $x_j$.
\label{prop:peripheral} \end{prop} 

\begin{proof} Let $c_i$ be an embedded loop in $R \cup \{x_i\}$ based at $x_i$ such that $c$ bounds an embedded disc $E$ in $S$ with
exactly one puncture $x_j$ of $R$ in its interior. Let $A$ be a closed annular neighborhood of $c$ in $S$ such that the boundary of $A$
is contained in $R$ and $x_i$ is the unique puncture of $R$ contained in the interior of $A$. Note that the boundary of $A$ consists of
two components, $C$ and $D$, where $D$ bounds an embedded disc $F$ contained in the interior of $E$, such that $x_j$ is the unique
puncture of $R$ in the interior of $F$, and $C$ bounds an embedded disc $G$ such that $x_i$ and $x_j$ are the unique punctures of $R$
contained in the interior of $G$. Note that $G = A \cup F$ and $A \cap F = D$. For simplicity, we may assume that the embedded loop $c$
is oriented so that $c$ is isotopic to the oriented circle $C$, where $C$ is given the boundary orientation induced from the annulus $A$
on the oriented surface $S$.  By the previous description of $\eta$, $\eta([c])$ is equal to the spin, $t_C \circ t_D^{-1}$ of $x_i$
about $c$. Since $D$ bounds a once punctured disc $F$ on $S$, $t_D$ is the trivial element of $PMod_R(S)$. Hence, $\eta([c]) = t_C$.
Since $C$ bounds the twice-punctured disc $G$ on $S$ and $x_i$ and $x_j$ are the unique punctures of $R$ in $G$, it follows that
$\eta([c])$ is the $2$-string twist braid $t_C$ with strings based at $x_i$ and $x_j$. 
\end{proof} 

Suppose that $x = x_i$. We shall denote $\mathcal{P}_i$ as $\mathcal{P}(x)$. 

\begin{prop} The restriction $\theta| : \pi_1(R \cup \{x\}, x) \rightarrow \pi_1(R \cup \{x\}, x)$ maps the peripheral subgroup
$\mathcal{P}(x)$ of $\pi_1(R \cup \{x\}, x)$ to $\mathcal{P}(x)$. 
\label{prop:presperipheral} \end{prop} 

\begin{proof} Let $f$ be an element of $\mathcal{P}(x)$. By Proposition \ref{prop:peripheral}, $f$ is a $2$-string twist braid
with strings based at $x$ and $u$, for some puncture $u$ of $R$ on $S$ with $u \neq x$. By Propositions \ref{prop:basic}, 
\ref{prop:bpbraids}, \ref{prop:twistbraids}, and \ref{prop:2twist}, the image of $f$ under $\theta$ is a $2$-string twist braid. On the other hand,
$\theta(f)$ is an element of $\pi_1(R \cup \{x\}, x)$. It follows, by an argument similar to the argument used in the proof of Proposition
\ref{prop:inducedisomorphism}, that the strings of $\theta(f)$ are based at $x$ and $v$, for some puncture $v$ of $R$ on $S$ with $v \neq x$.
Hence, by Proposition \ref{prop:peripheral},
$\theta(f)$ is an element of $\mathcal{P}(x)$. 

This proves that $\theta| : \pi_1(R \cup \{x\}, x) \rightarrow \pi_1(R \cup \{x\}, x)$ maps $\mathcal{P}(x)$ into $\mathcal{P}(x)$.
Likewise, the restriction $\theta^{-1}| : \pi_1(R \cup \{x\}, x) \rightarrow \pi_1(R \cup \{x\}, x)$ maps $\mathcal{P}(x)$ into 
$\mathcal{P}(x)$. 

Since $\theta^{-1}| : \pi_1(R \cup \{x\}, x) \rightarrow \pi_1(R \cup \{x\}, x)$ is the inverse of 
$\theta| : \pi_1(R \cup \{x\}, x) \rightarrow \pi_1(R \cup \{x\}, x)$, the result follows. 
\end{proof}

\begin{prop} There exists a homeomorphism $\overline{G}: (R \cup \{x\}, x) \rightarrow (R \cup \{x\}, x)$ such that the automorphism
$\theta| : \pi_1(R \cup \{x\}, x) \rightarrow \pi_1(R \cup \{x\}, x)$ of $\pi_1(R \cup \{x\}, x)$ is induced by $\overline{G}$. 
\label{prop:indautgeom} \end{prop} 

\begin{proof} By Propositions \ref{prop:inducedautomorphism} and \ref{prop:presperipheral}, $\theta| : \pi_1(R \cup \{x\}, x)
\rightarrow \pi_1(R \cup \{x\}, x)$ is an automorphism of $\pi_1(R \cup \{x\}, x)$ preserving the peripheral structure $\mathcal{P}(x)$
of $\pi_1(R \cup \{x\}, x)$. By a well-known result about automorphisms of fundamental groups of surfaces, it follows that $\theta|$
is induced by the restriction $G : R \rightarrow R$ of a homeomorphism $\overline{G} : (R \cup \{x\}, x) \rightarrow (R \cup \{x\}, x)$
of the pointed space $(R \cup \{x\}, x)$ (cf. \cite{zvc}, Theorems 3.3.11, 5.7.1, and 5.13.1).
\end{proof} 

\begin{prop} There exists a homeomorphism $\overline{G}: (R \cup \{x\}, x) \rightarrow (R \cup \{x\}, x)$ such that the
automorphism $\theta| : \pi_1(R \cup \{x\}, x) \rightarrow \pi_1(R \cup \{x\}, x)$ of $\pi_1(R \cup \{x\}, x)$ is induced by
$\overline{G}$. 
\label{prop:indautgeom} \end{prop} 

\begin{proof} By Propositions \ref{prop:inducedautomorphism} and \ref{prop:presperipheral}, $\theta| : \pi_1(R \cup \{x\}, x)
\rightarrow \pi_1(R \cup \{x\}, x)$ is an automorphism of $\pi_1(R \cup \{x\}, x)$ preserving the peripheral structure $\mathcal{P}(x)$
of $\pi_1(R \cup \{x\}, x)$. By a well-known result about automorphisms of fundamental groups of surfaces, it follows that $\theta|$
is induced by the restriction $G : R \rightarrow R$ of a homeomorphism $\overline{G} : (R \cup \{x\}, x) \rightarrow (R \cup \{x\}, x)$
of the pointed space $(R \cup \{x\}, x)$ (cf. \cite{zvc}, Theorems 3.3.11, 5.7.1, and 5.13.1).
\end{proof} 

\begin{prop} Let $\overline{G}: (R \cup \{x\}, x) \rightarrow (R \cup \{x\}, x)$ be as in Proposition \ref{prop:indautgeom}.
Let $G : R \rightarrow R$ be the restriction of $\overline{G}$ to $R$. Let $G_* : PMod_R(S) \rightarrow PMod_R(S)$ be the automorphism
of $PMod_R(S)$ defined by the rule $G_*([H]) = [G \circ H \circ G^{-1}]$ for each mapping class $[H]$ in $PMod_R(S)$. Let $\phi :
PMod_R(S) \rightarrow PMod_R(S)$ be the composition $(G_*)^{-1} \circ \theta : PMod_R(S) \rightarrow PMod_R(S)$. Then the restriction of
$\phi$ to $\pi_1(R \cup \{x\}, x)$ is equal to the identity. 
\label{prop:trivonsurfsubgroup} \end{prop} 

\begin{proof} This is immediate from Propositions \ref{prop:equivariance} and \ref{prop:indautgeom}. 
\end{proof}

\section{the main results} 
\label{sec:mainresults} 

\begin{thm} Let $S$ be a closed connected orientable surface of positive genus. Let $n$ be an integer greater than $2$.
Let $R$ be a surface obtained from $S$ by removing $n$ distinct points from $S$. Let $\chi : PMod_R(S) \rightarrow PMod_R(S)$ be an
automorphism of the pure braid group $PMod_R(S)$. Then there exists a homeomorphism $K: R \rightarrow R$ of $R$ such that $\chi : PMod_R(S)
\rightarrow PMod_R(S)$ is equal to the automorphism $K_* : PMod_R(S) \rightarrow PMod_R(S)$ which is defined by the rule $K_*([H]) = [K
\circ H \circ K^{-1}]$ for each mapping class $[H]$ in $PMod_R(S)$.  
\label{thm:purebraidthm} \end{thm} 

\begin{proof} Let $x$ be a puncture of $R$ on $S$. As in Proposition \ref{prop:inducedautomorphism}, we
may choose a homeomorphism $F : R \rightarrow R$ of $R$ such that the composition $\theta$ of $F_*$ and $\chi$ restricts to an
automorphism $\theta|$ of $\pi_1(R \cup \{x\}, x)$. By Propositions \ref{prop:indautgeom} and \ref{prop:trivonsurfsubgroup}, we may
choose a homeomorphism $G : R \rightarrow R$ of $R$ such that the composition $\phi$ of $(G_*)^{-1}$ and $\theta$ fixes each element of
$\pi_1(R \cup \{x\}, x)$. It follows, by Proposition \ref{prop:surfsubgp}, that $\phi : PMod_R(S) \rightarrow PMod_R(S)$
is equal to the identity automorphism $id : PMod_R(S) \rightarrow PMod_R(S)$ of $PMod_R(S)$. Hence, 

$$id = \phi = (G_*)^{-1} \circ \theta = (G_*)^{-1} \circ (F_*) \circ \chi$$. 

Hence, $\chi = (F^{-1} \circ G)_*$. Let $K = F^{-1} \circ G$. Then $K : R \rightarrow R$ is a homeomorphism of $R$ and $\chi = K_*$.
This completes the proof of Theorem \ref{thm:purebraidthm}.
\end{proof}

Our main result for the braid group $Mod_R(S)$ will follow from Theorem \ref{thm:purebraidthm} by using Theorem $2$ of \cite{i3} and 
the following proposition. 

\begin{prop} Suppose that $\tau : Mod_R(S) \rightarrow Mod_R(S)$ is an automorphism of $Mod_R(S)$ such that $\tau$ fixes each element of
$PMod_R(S)$. Then $\tau$ is the identity automorphism of $Mod_R(S)$. 
\label{prop:triviftrivonpurebraids} \end{prop} 

\begin{proof} Suppose that $x$ and $y$ are distinct punctures of $R$ on $S$. Let $D$ be an embedded disc on $S$ such that the
boundary of $D$ is contained in $R$ and $x$ and $y$ are the unique punctures of $R$ on $S$ which are contained in the interior of
$D$. Choose a homeomorphism $F : R \rightarrow R$ which is supported on $D$ and interchanges $x$ and $y$ in such a way that $f^2$ is
equal to $t_C$, where $C$ is the boundary of $D$ and $f$ is the mapping class of $F$. Since $F$ is supported on $D$, $f$ is an element
of $PMod_R(S)$. Following common terminology, we shall call $f$ an elementary braid. Recall that $Mod_R(S)$ is generated by $PMod_R(S)$ and
elementary braids. Hence, it remains only to show that $\tau$ fixes each elementary braid, $f$. 

Let $C_i$, $1 < i \leq n$ be a sequence of disjoint circles in $R$ such that $C_2 = C$ and $C_i$ bounds an embedded disc $D_i$ on $S$ such
that the interior of $D_i$ in $S$ contains exactly $i$ punctures of $R$ on $S$. Note that $D_2 = D$. Let $t_i = t_{C_i}$. Then $t_i$ is an
$i$-string twist braid on $R$ supported on $D_i$ and, hence, an element of $PMod_R(S)$.  

Note that $f^2 = t_C = t_{C_2} = t_2$. Hence, $f$ commutes with $t_2$. Suppose that $2 < i \leq n$. Note that $D_{i - 1}$ is contained in the interior of
$D_i$. Hence, $C_i$ is contained in the complement of the support $D$ of $F$. It follows that $f$ commutes with $t_i$. 

Let $2 \leq i \leq n$. Then $t_i$ is an element of $PMod_R(S)$ and $f$ commutes with $t_i$. Hence, $\tau(t_i) = t_i$ and $\tau(f)$ commutes with
$\tau(t_i)$. This implies that $\tau(f)$ commutes with $t_i$. It follows that $\tau(f)$ is the mapping class $[P]$ of a homeomorphism $P : R
\rightarrow R$ such that $P(C_i) = C_i$ for all $i$ with $2 \leq i \leq n$. Note that $f$ is the mapping class $F$ of a homeomorphism $F : R
\rightarrow R$ such that $F(C_i) = C_i$ for all $i$ with $2 \leq i \leq n$.

Let $p_i$, $1 \leq i \leq n$, denote the punctures of $R$ on $S$. We may assume that $p_1 = x$, $p_2 = y$, and $p_i$ is contained in the interior of
$D_i$ for $2 \leq i \leq n$. We may extend $P$ to a homeomorphism $\overline{P} : S \rightarrow S$ of $S$. Note that $\overline{P}(p_i) = p_i$ for all
$i$ with $2 < i \leq n$. Hence, either $\overline{P}(x) = x$ or $\overline{P}(x) = y$. Suppose that $\overline{P}(x) = x$. Then $\overline{P}(y) = y$.
This implies that $\tau(f)$, the mapping class of $P$, is an element of $PMod_R(S)$. Since $\tau$ fixes each element of $PMod_R(S)$, this implies that
$f$ is an element of $PMod_R(S)$. On the other hand, since $f$ is an elementary braid ``switching'' $x$ and $y$, $f$ is not an element of $PMod_R(S)$.
This is a contradiction. Hence, $\overline{P}(x) = y$. This implies that $\overline{P}(y) = x$. This proves that the actions of $f$ and $\tau(f)$ on the
set of punctures of $R$ on $S$ are equal. 

Let $g = f^{-1} \tau(f)$. Note that $g$ is the mapping class of the homeomorphism $G : R \rightarrow R$, where $G = F^{-1} \circ P$. Note that $G(C_i) =
C_i$ for all $i$ with $2 \leq i \leq n$. Moreover, $G$ extends to a homeomorphism $\overline{G} : S \rightarrow S$ such that $G(p_i) = p_i$ for all $i$
with $2 \leq i \leq n$. Hence, $g$ is an element of $PMod_R(S)$ fixing the isotopy classes $[C_i]$, $1 \leq i \leq n$. It follows that
$g$ is the composition of a multitwist $t$ on the circles $C_i$, $1 \leq i \leq n$ with a mapping class $h$ supported on the complement $K$ of the
interior of $D_n$. Note that each of the mapping classes $t$ and $h$ are supported on the complement of the interior of $D$. Hence, $g$ commutes with
$f$. Since $g = f^{-1} \tau(f)$, it follows that $f$ commutes with $\tau(f)$. 

We conclude that $g^2 = f^{-1} \tau(f) f^{-1} \tau(f) = f^{-2} (\tau(f))^2 = t_2^{-1} (\tau(f))^2$. Since $f^2 = t_2$ $(\tau(f))^2  =
\tau(f^2) = \tau(t_2) = t_2$. Hence, $g^2 = id$. Since $g$ is a pure braid, it follows by Theorem \ref{thm:pureispure}, that $g = id$. That is to say,
$\tau(f) = f$.  
\end{proof}

\begin{thm} Let $S$ be a closed connected orientable surface of genus greater than $1$. Let $n$ be an integer greater than $2$.
Let $R$ be a surface obtained from $S$ by removing $n$ distinct points from $S$. Let $\chi : Mod_R(S) \rightarrow Mod_R(S)$ be an automorphism 
of the braid group $Mod_R(S)$. Then there exists a homeomorphism $F: R \rightarrow R$ of $R$ such that $\chi : Mod_R(S) \rightarrow
Mod_R(S)$ is equal to the automorphism $F_* : Mod_R(S) \rightarrow Mod_R(S)$ which is defined by the rule $F_*([H]) = [F \circ H \circ
F^{-1}]$ for each mapping class $[H]$ in $Mod_R(S)$. 
\label{thm:braidthm} \end{thm}

\begin{proof} By Theorem $2$ of \cite{i3}, $PMod_R(S)$ is a characteristic subgroup of $Mod_R(S)$. Hence, $\chi : Mod_R(S) \rightarrow
Mod_R(S)$ restricts to an automorphism $\chi| : PMod_R(S) \rightarrow PMod_R(S)$. By Theorem \ref{thm:purebraidthm}, there exists a
homeomorphism $K : R \rightarrow R$ of $R$ such that $\chi| : PMod_R(S) \rightarrow PMod_R(S)$ is equal to $K_* : PMod_R(S) \rightarrow
PMod_R(S)$. We shall show that $\chi : Mod_R(S) \rightarrow Mod_R(S)$ is equal to $K_* : Mod_R(S) \rightarrow Mod_R(S)$. 

Let $\tau : Mod_R(S) \rightarrow Mod_R(S)$ be the composition $K_*^{_1} \circ \chi : Mod_R(S) \rightarrow Mod_R(S)$. Since 
$\chi| : PMod_R(S) \rightarrow PMod_R(S)$ is equal to $K_* : PMod_R(S) \rightarrow PMod_R(S)$, $\tau$ fixes each element of $PMod_R(S)$. Hence, by 
Proposition \ref{prop:triviftrivonpurebraids}, $\tau$ is the identity automorphism of $Mod_R(S)$. That is to say, $\tau : Mod_R(S) \rightarrow Mod_R(S)$
is equal to $K_* : Mod_R(S) \rightarrow Mod_R(S)$. 
\end{proof}


\begin{thebibliography}{XXXX}

\bibitem[B]{b} J. S. Birman, {\it Braids, Links and Mapping Class Groups}, Annals of Math. 
Studies, No. 82, Princeton University Press, Princeton, New Jersey, 1974.

\bibitem[BH]{bh} J. S. Birman and  H. M. Hilden, On isotopies of homeomorphisms of 
Riemann surfaces, Ann. of Math. Vol. 97, No. 3, (1973), 424-439.

\bibitem[BLM]{blm} J. S. Birman, A. Lubotzky, and J. McCarthy, Abelian and solvable subgroups 
of the mapping class group, Duke Math. J. V. 50 (1983), 1107-1120.

\bibitem[D]{d} M. Dehn, Die Gruppe der Abbildungsklassen, Acta Math. V. 69 (1938), 135 - 206.

\bibitem[DG]{dg} J. L. Dyer and E. K. Grossman, {\it The automorphism groups of the braid groups}, Amer. J. Math. 103 (1981), 1151-1169.

\bibitem[EE]{ee} C. J. Earle and J. Eells, {\it A fibre bundle description of Teichm\"{u}ller theory}, J. Diff. 
Geometry, V. 3, No. 1 (1969), 19-43.

\bibitem[F]{f} B. Farb, {it Automorphisms of the Torelli group}, to appear, preliminary report - AMS Sectional Meeting, Ann Arbor, 
Michigan, March 1, 2003

\bibitem[FLP]{flp} A. Fathi, F. Laudenbach. and V. Po\'{e}naru, {\it Travaux de Thurston 
sur les surfaces}, S\'eminaire Orsay, Ast\'{e}risque, Vol. 66-67, Soc. Math. de France, 1979.

\bibitem[GS]{gs} S. M. Gersten, J. R. Stallings, {\it Combinatorial Group Theory and Topology}, 
Annals of Math. Studies No. 111, Princeton University Press, Princeton, New Jersey, 1987.

\bibitem[GG]{gr} M. Gromov, Hyperbolic groups, in: S. M. Gersten, Editor, {\it Essays in group theory}, 
MSRI Publications, Vol. 8, Springer-Verlag, Berlin, Heidelberg, New York, 1987, pp. 75-263.

\bibitem[G]{grossman} E. Grossman, On the residual finiteness of certain mapping class groups,
J. London Math. Soc. V. 9 (1974), 160-164.

\bibitem[I1]{i1} N. V. Ivanov, Automorphisms of Teichm\"{u}ller modular groups, Lecture Notes 
in Math., No. 1346, Springer-Verlag, Berlin and New York, 1988, 199-270.

\bibitem[I2]{i2} N. V. Ivanov, {\it Subgroups of Teichm\"{u}ller Modular Groups}, 
Translations of Mathematical Monographs, Vol. 115, American Math. Soc., Providence, Rhode Island, 
1992.

\bibitem[I3]{i3} N. V. Ivanov, {\it Permutation representations of braid groups of surfaces.} (Russian) Mat. Sb. 181 (1990), no. 11,
1464--1474; translation in Math. USSR-Sb. 71 (1992), no. 2, 309--318

\bibitem[IM]{im} N. V. Ivanov and J. D. McCarthy, {\it On injective homomorphisms between Teichm\"{u}ller modular groups, I}, 
Invent. math. 135, 425-486 (1999)

\bibitem[J]{j} D. L. Johnson, Homeomorphisms of a surface which act trivially on homology, 
Proc. Amer. Math. Soc. V. 75 (1979), 119-125.

\bibitem[K]{k} M. Korkmaz, {\it Automorphisms of complexes of curves on punctured spheres and on punctured tori}, Topology Appl. 95 (1999),
no. 2, 85--111

\bibitem[M]{m} J. D. McCarthy, Automorphisms of surface mapping class groups. A recent 
theorem of N. Ivanov, Invent. Math. V. 84, F. 1 (1986), 49-71.

\bibitem[MV]{mv} J. D. McCarthy and W. R. Vautaw, {\it Automorphisms of the Torelli Group in genus $>$ 2}, to appear, preliminary report - AMS
Sectional Meeting, Baton Rouge, Louisiana, March 14, 2003

\bibitem[ZVC]{zvc} H. Zieschang, E. Vogt, H.-D. Coldeway, {\it Surfaces and planar discontinuous groups}, Lecture Notes 
in Math., No. 835, Springer-Verlag, Berlin, Heidelberg, New York, 1980. 

\end{thebibliography}
\end{document}